\documentclass{article}
\usepackage{amsmath,amssymb}
\usepackage[all]{xy}
\font\ibf=cmbxti10

\title{Discrete Symmetry with Compact Fundamental Domain, and Geometric Simple Connectivity \\
-- {\large A Provisional Outline of Work in Progress} -- }
\author{Valentin {\sc Po\'enaru}\footnote{Professor Emeritus at the Universit\'e Paris Sud-Orsay, Math\'ematiques 425, 91405 Orsay Cedex, France. e-mail: valpoe@hotmail.com}}
\date{October 2007}

\begin{document}

\maketitle

\vglue 1cm

\begin{abstract}
We show that a certain geometric property, the QSF introduced by S.~Brick and M.~Mihalik, is universally true for {\ibf all} finitely presented groups $\Gamma$. One way of defining this property is the existence of a smooth compact manifold $M$ with $\pi_1 \, M = \Gamma$, such that $\tilde M$ is geometrically simply-connected ({\it i.e.} without handles of index $\lambda = 1$). There are also alternative, more group-theoretical definitions, which are presentation independent. But $\Gamma \in {\rm QSF}$ is not only a universal property, it is quite highly non-trivial too; its very special case for $\Gamma = \pi_1 \, M^3$ (where it means $\pi_1^{\infty} \, \tilde M^3 = 0$) is actually already known, as a corollary of G.~Perelman's big breakthrough on the Geometrization of 3-Manifolds.
\end{abstract}

\vglue 1cm

\section{Introduction}

This is a paper about {\ibf all} finitely presented groups, and no others will ever be considered here. Our main result is the following.

\bigskip

\noindent {\bf Theorem 1.} {\it For {\ibf any} (finitely presented) group $\Gamma$, there is a finite simplicial complex $P$ with $\pi_1 \, P = \Gamma$, such that the universal covering space $\tilde P$ has the following property}

\medskip

\noindent (QSF) {\it For any compact $k \subset \tilde P$ there is an (abstract) simply-connected compact $K$, coming with an inclusion $k \overset{i}{\longrightarrow} K$ and with a continuous map $F$, inducing a commutative diagram}
$$
\xymatrix{
k \ar[rr]^{i} \ar[dr]_{\underset{\mbox{\scriptsize inclusion}}{\rm the \, natural}}  &&K \ar[dl]^{F}  \\ &\tilde P
}
$$
{\it such that} $\{$the double point set $M_2 (F) \subset K \} \cap i(k) = \emptyset$. {\it It is understood here that the diagram above is meant in the simplicial category.}

\bigskip

The QSF concept (meaning ``quasi-simply-filtered'') has been introduced by S.~Brick and M.~Mihalik \cite{B-M}; but see here J.~Stallings \cite{S} too. It is a notion which makes sense, {\it a priori}, for locally compact spaces $X$. But then, it is shown in \cite{B-M} that if $P_1 , P_2$ are two compact spaces with $\pi_1 = \Gamma$, then $\tilde P_1 \in {\rm QSF} \Longleftrightarrow \tilde P_2 \in {\rm QSF}$. This means that Theorem~1 is satisfied for {\it any} $P$ such that $P$ is compact and $\pi_1 \, P = \Gamma$. It also means that the (QSF) is a {\it bona fide} group theoretical notion, which is presentation independent; most of the other naturally related notions are not. Our Theorem~1 answers affirmatively a question raised at the very end of \cite{B-M}. The reader should probably also know that the classical Whitehead 3-manifold ${\rm Wh}^3$ is a locally compact space which, indeed, is not QSF.

\medskip

As far as I am concerned, there are at least two good reasons which make that our theorem above is not without some interest. To begin with, on a more philosophical side, here we have a property which is true for all groups, at least so our theorem claims, but then as I shall soon show, this property is, also, very highly non-trivial. That such a combination of features should be possible, at all, goes against the standard wisdom in the field. So, how non-trivial is, actually, our theorem?

\medskip

Let us consider the very special case when $\Gamma = \pi_1 \, M^3$, where $M^3$ is a closed 3-manifold. It is shown in \cite{B-M} that, in this special case $\pi_1 \, M^3 \in {\rm QSF}$ iff $\pi_1^{\infty} \, \tilde M^3 = 0$ (which, in purely group-theoretical terms can be also formulated as $\pi_1^{\infty} \, \Gamma = 0$). But then, the ``$\pi_1^{\infty} \, \tilde M^3 = 0$'' also has a purely topological interpretation and the off-shot now is that, {\it via} \cite{B-M}, the following statement is now a {\it consequence}, or corollary, of Theorem~1, namely
$$
\mbox{For any closed irreducible $M^3$ with infinite $\pi_1$, we have $\tilde M^3 = R^3$.} \leqno (1)
$$
\noindent Now, of course, this (1), which for long time has been a celebrated conjecture, is today already known, since it also follows from the Thurston Geometrization of 3-Manifolds, {\it i.e.} it is implied by Grisha Perelman's big breakthrough \cite{Pe1}, \cite{Pe2}, \cite{Pe3}, \cite{Mo}, \cite{Be}, \cite{Bes}, \cite{Ma}.

\medskip

The point which I am trying to make here is that the {\it only} known way to get the (1) above, other than applying Theorem~1, is to use Perelman's spectacular work. This certainly also shows that the fact stated by Theorem~1 is quite highly non-trivial. 

\medskip

As a side-remark, it is only for groups $\Gamma$ which are of the form $\Gamma = \pi_1 \, M^3$, that $\pi_1^{\infty} = 0$ and QSF are equivalent properties. Before closing this first, philosophical comment concerning our theorem let me also add that I believe that there should be a way to reconcile our Theorem~1 with the standard wisdom. The question raised at the very end of Section 3 below, may be on the right track, as far as this issue is concerned; after all, when ``discrete symmetries'' at large are concerned, then there is more under the Sun than just ``all the groups''.

\medskip

But then, even closer to my own circle of ideas (see here, for instance, \cite{Po7}, \cite{Po8}, \cite{Po9}, \cite{Po10}), there is also a second reason why, for me at least, Theorem~1 is of interest.

\medskip

This is because it reveals a relatively surprising connection between discrete groups theory and the concept of {\ibf geometric simple connectivity}, which comes from a very far away mathematical universe. I remind the reader that, by definition, a smooth manifold is GSC ($=$ Geometrically Simply Connected) if it admits a smooth handle-body decomposition without handles of index $\lambda = 1$. A very pedestrian introduction to this topic is to be found in \cite{Po-Ta1}. We will also soon see that the GSC concept makes sense in more general contexts.

\medskip

Anyway, at a first, more simple-minded level, here is

\bigskip

\noindent {\bf Corollary 2.} {\it For any group $\Gamma$ and for sufficiently large $n \in Z_+$, there is a smooth compact $n$-manifold $M^n$ such that $\pi_1 \, M^n = \Gamma$ and $\tilde M^n \in {\rm GSC}$.}

\bigskip

This Corollary~2 follows by putting together Theorem~1 and results of L.~Funar and D.~Otera \cite{Fu-O}, \cite{O} (and see the \cite{Fu-G} too). The reader should notice here that, when groups are concerned, while QSF is a presentation-independent notion, GSC is not. Since GSC immediately implies QSF, one may view the work of Funar and Otera as showing that, in an appropriately defined set-up, there is equivalence between ``$\Gamma \in {\rm QSF}$'' and ``$\Gamma \in {\rm GSC}$''.

\medskip

At a deeper level, the proof of Theorem~1 relies very heavily on a certain MAIN LEMMA, specifically this is Lemma~6 below, and what this lemma says, is that a certain smooth, very high-dimensional non-compact manifold with large non-empty boundary, very much connected with $\tilde M^3 (\Gamma)$ (which is a particular finite presentation of $\Gamma$, to be soon introduced), but without being exactly a high-dimensional thickening of it, is GSC. We will have a lot to say later concerning this manifold, called $S_u (\tilde M^3 (\Gamma))$. For right now, it suffices to say that it comes naturally equipped with a free action of $\Gamma$ but for which, unfortunately so to say, the fundamental domain $S_u (\tilde M^3 (\Gamma))/\Gamma$ fails to be compact. If one thinks of the $\tilde P$ which occurs in our theorem as being made out of fundamental domains which are like solid compact ice-cubes, then the $S_u (\tilde M^3 (\Gamma))$ itself is gotten, essentially, by replacing each of these compact ice-cubes by a non-compact infinitely foamy structure, and then going to a very high-dimensional thickening of this foam.

\medskip

In order to clinch the proof of our Theorem~1, one finally also needs to prove the following kind of implication

$$
\{ S_u (\tilde M^3 (\Gamma)) \in {\rm GSC}, \ \mbox{which is what the MAIN LEMMA claims}\} \Longrightarrow \leqno (2)
$$
$$
\Longrightarrow \{ \Gamma \in {\rm QSF}\} \, .
$$

At the level of this short outline there is no room even for a sketch of proof of the implication (2) above. Suffices to say here that this (2) is a considerably easier step than the main lemma itself, about which a lot will be said below.

\vglue 2mm

$$
* \ * \ * \ * \ *
$$

\vglue 3mm

The many criticisms, comments and suggestions which David Gabai has made in connection with my earlier ill-fated attempt of proving $\pi_1^{\infty} \, \tilde M^3 = 0$ (circa 2000), were essential for the present work, which could not have existed without them. Like in other occasions too, his help was crucial for me.

\medskip

Thanks are also due to Louis Funar and Daniele Otera for very useful conversations. Actually, it was Daniele who, at that time my PhD student, first told me about QSF, and who also insisted that I should look into it.

\medskip

The friendly help of the IH\'ES made the existence of the present paper possible, and I wish to thank here C\'ecile Cheikhchoukh for the typing.

\section{GSC Pseudo-spine representations of groups}

For technical reasons which will become clear later on, when flesh and bone are put onto $\Gamma$, by picking up a compact $P$ with $\pi_1 \, P = \Gamma$, we will make the following kind of choice for our $P$. We will disguise our general $\Gamma$ as a 3-manifold group, in the sense that we will work with a compact bounded {\ibf singular} 3-manifold $M^3 (\Gamma)$, with $\pi_1 \, M^3 (\Gamma) = \Gamma$. Typically, such an $M^3 (\Gamma)$ will be a {\it bona fide} finite handle-body $H^3$ with finitely many 2-handles attached. But the map from the attaching zones $\sum S^1 \times I$ into $\partial H^3$ will be allowed to have generic double points. The corresponding little squares (contained inside $\partial H^3$) are called the {\ibf immortal singularities} ${\rm Sing} \, (M^3 (\Gamma)) \subset M^3 (\Gamma)$ of our $M^3 (\Gamma)$.

\medskip

Our singular $M^3 (\Gamma)$'s in question have absolutely no vocation for getting any kind of ``geometric structures'', beyond the very austere one given by the free action of $\Gamma$ on $\tilde M^3 (\Gamma)$, in other words other than $\Gamma$ itself.

\medskip

The technology which was first introduced in \cite{Po1} and then used in subsequent papers like \cite{Po2}, \cite{Po3}, and others, although it was very opportunistically formulated there in terms of smooth 3-manifolds, is very easily adaptable to a vast variety of more general situations. In particular, it can be adapted to our present $M^3 (\Gamma)$'s and their univeral covering spaces $\tilde M^3 (\Gamma)$. For the convenience of the reader, the basic necessary material will be very briefly reviewed here.

\medskip

One considers now infinite, not necessarily locally finite complexes $X$, of dimension $\leq 3$, endowed with non-degenerate cellular maps $f$ into $M^3 (\Gamma)$ or $\tilde M^3 (\Gamma)$; non-degenerate means here that $f$ injects on the individual cells. The points of $X$ where $f$ is not immersive, are by definition the {\ibf mortal singularities} ${\rm Sing} \, (f) \subset X$.

\medskip

Two kinds of equivalence relations will be associated to such non-degenerate maps $f$, namely the
$$
\Psi (f) \subset \Phi (f) \subset X \times X \, . \leqno (3)
$$
Here $\Phi (f)$ is the simple-minded equivalence relation where $(x,y) \in \Phi (f)$ means that $f(x) = f(y)$. The more subtle and not so easily definable $\Psi (f)$, is the smallest equivalence relation compatible with $f$, which kills all the mortal singularities. It can be shown that an unambiguous mathematical meaning can be given to such a definition, that it comes with a commutative diagram
$$
\xymatrix{
X \ar[rr]^{f} \ar[dr]_{\pi (\Psi (f))}  &&&\!\!\!\!\!\!\!\!\!\!\!\!\!\!\!\{M^3 (\Gamma) \ {\rm or} \ \tilde M^3 (\Gamma), \mbox{according to the case} \}  \\ &X/\Psi (f) \ar[ur]_{f_1 (\Psi (f))}
} \leqno (4)
$$
where $f_1$ is now an immersion ($\Leftrightarrow {\rm Sing} \, (f_1) = \emptyset$), and also that the following map is a surjection
$$
\xymatrix{
\pi_1 \, X \ar[rr]^{\!\!\!\!\!\!\!\!\!\!\pi_*}  && \ \pi_1 (X / \Psi (f)) \, .   \\ &
} \leqno (5)
$$
\vglue -5mm
\noindent Moreover, the map $X \overset{\pi}{\longrightarrow} X/\Psi (f)$ is realizable (not uniquely) {\it via} a sequence of folding maps and any such sequence will be called a {\bf zipping} (or a zipping strategy). For the time being $X$ is 3-dimensional.

\bigskip

\noindent {\bf Lemma 3.} {\it There is an $X$ which is {\rm GSC}, coming with a non degenerate surjection $X \overset{f}{\longrightarrow} \tilde M^3 (\Gamma)$, which is such that one should have}
$$
\Psi (f) = \Phi (f) \, . \leqno (6)
$$

\bigskip

Before anything else, here are some comments concerning this lemma.

\medskip

A) Without anything else being required beyond GSC and $\Psi = \Phi$, any pair $(X,f)$ like in our lemma, will be called a {\ibf representation} for $\Gamma$, or for $\tilde M^3 (\Gamma)$. These are the representations meant in the title of this section. Several, better and better representations will be introduced below. It so happens that this initial one, offered by the lemma above comes with an $f$ which is surjective. Also, although not locally finite our $X$ from Lemma~3 is relatively simple-minded. For our next representations, the source will be more exotic and the map will fail to surject. In Theorem~4, for instance, the target space $\tilde M^3 (\Gamma)$ will only be the closure of the image of the map in question; later on, things will get even worse, from this viewpoint.

\medskip

B) In the context of Lemma~3 itself, but not later on, we can actually do better than just GSC. The $X$ in that lemma is actually arborescent, {\it i.e.} gettable from a point {\it via} an infinite sequence of Whitehead dilatation (of dimension $\leq 3$), and this implies GSC.

\medskip

On the other hand, notice that we are using now the GSC concept outside of the smooth context, hence some explanations may be necessary. So, we consider now a context, either of infinite cell-complexes or of smooth non-compact manifolds, with possibly $\partial \ne \emptyset$. We will say that our corresponding object, call it still $X$, is GSC if it has a cell-decomposition or handle-body decomposition, according to the case, of the following type. Start with a PROPERLY embedded tree (or with a smooth regular neighbourhood of such a tree); with this should come now a cell-decomposition, or handle-decomposition
$$
X = T \cup \sum_{1}^{\infty} \{ \mbox{$1$-handles (or $1$-cells)} \, H_i^1 \} \cup \sum_{1}^{\infty} H_j^2 \cup \biggl\{ \sum_{k,\lambda \geq 2} H_k^{\lambda} \biggl\} \, , \leqno (7)
$$ 
such that the {\ibf geometric intersection matrix}, which counts without any kind of $\pm$ signs, how many times $H_j^2$ touches $H_i^1$, takes the following ``easy'' id $+$ nilpotent form
$$
H_j^2 \cdot H_i^1 = \delta_{ji} + a_{ji} \, , \quad \mbox{where} \quad a_{ji} \in Z_+ \quad \mbox{and} \quad a_{ji} > 0 \Rightarrow j > i \, . \leqno (8)
$$

\medskip

C) (additional comments to B)) In this paper we will distinguish between PROPER, meaning inverse image of compact is compact, and proper meaning (inverse image of boundary) $=$ boundary.

\medskip

Also, there is a notion of ``difficult'' id $+$ nilpotent, gotten by reversing the last inequality in (8), and this is no longer GSC. For instance, the Whitehead manifold ${\rm Wh}^3$, several times mentioned in this paper, always impersonating one of the villains of the cast, admits handle-body decompositions of the difficult id $+$ nilpotent type. But then, ${\rm Wh}^3$ is certainly not GSC either.

\medskip

D) The $X$ in any of the representations which we are considering,  has only mortal singularities, just like the target $\tilde M^3 (\Gamma)$ only has immortal ones. This latter ones are created by the zipping process itself, {\it i.e.} by the following arrow, which is always to be thought of as coming with a precise factorization into successive folding maps
$$
\xymatrix{
X \ar[rr]^{\!\!\!\!\!\!\!\!\!\!\!\!\!\!\!\!\!\!\!\!\!\!\!\pi}   &&X / \Psi (f) = \tilde M^3 (\Gamma) \, .   \\ &
} \leqno (9)
$$
\vglue -5mm
At any of its infinitely many singularities, the $X$ from Lemma~3 looks, locally like the $\{$Riemann surface of $\log z \} \times R$.

\medskip

After these preliminaries, we move now to more serious items.

\bigskip

\noindent {\bf Theorem 4.} {\it There is a $3^{\rm d}$ representation
$$
\xymatrix{
Y(\infty) \ar[rr]^{G(\infty)}   &&\tilde M^3 (\Gamma) \, ,   \\ &
} \leqno (10)
$$
\vglue -5mm
\noindent which also has the following additional features.}

\medskip

1) {\it The source $Y(\infty)$ is now locally finite; we will call this our first finiteness condition.}

\medskip

2) {\it The representation $(10)$ is EQUIVARIANT. Specifically, the representation space $Y(\infty)$ itself comes equipped with a free action of $\Gamma$ and then, for each $x \in Y(\infty)$ and $g \in \Gamma$, we have also that}
$$
G(\infty) (gx) = g \, G (\infty) (x) \, .
$$

\medskip

3) {\it There is now a {\ibf uniform bound} $M > 0$ and also a zipping strategy for $(10)$, like in $(9)$ with the appropriate change of notations, such that for any}
$$
(x,y) \in \Psi (G(\infty)) = \Phi (G(\infty)) \subset Y(\infty) \times Y(\infty)
$$
{\it there is a {\ibf zipping path} $\lambda (x,y) \subset \Psi (G(\infty)) \cup {\rm Diag} \, ({\rm Sing} \, (G(\infty)) \subset Y(\infty) \times Y(\infty))$, connecting our $(x,y)$ to the mortal singularities of $G(\infty)$, such that}
$$
{\rm length} \, \lambda (x,y) \leq M \, . \leqno (11)
$$

\bigskip

Here are some comments  concerning this statement. To begin with, it turns out that our previous paper \cite{Po-Ta2}, although written originally for smooth 3-manifolds, is perfectly well adaptable for the context of the singular $M^3 (\Gamma)$. It is actually like tailor-made for the 1) $+$ 2) above. Just like in \cite{Po-Ta2}, our $Y(\infty)$ has a sort of {\ibf train-track} handle-body decomposition put together, very singularly, out of 3-dimensional handles of index $\lambda \leq 2$, deprived each of its lateral surface ($=$ boundary minus attaching zone), and hence rendered non-compact.

\medskip

It should not be hard for the reader to figure out what ``zipping paths'' means and then, there are various natural options for measuring their lengths. One may count, for instance, the number of necessary folding operations. Or, one may also lift some equivariant riemannian-type metric from $\tilde M^3 (\Gamma)$ to $\Psi (G(\infty)) \subset Y(\infty) \times Y(\infty)$, and then use it for measuring lengths. All these various notions of zipping length are equivalent up to quasi-isometry.

\medskip

All these things having been said, the proof of the 3) requires a new, infinitistic argument. We start with a representation $(Y(1), G(1))$ having only 1) $+$ 2) but not necessarily 3). Then, without loosing 1) + 2) we try to improve things, {\it via} successive enlargements. One goes first to an enlargement $(Y(2), G(2))$, for which the zipping lengths are uniformly bounded when restricted to $\Psi (G(1)) = \Phi (G(1)) \subset \Psi (G(2))$. Next, one enlarges again to $(Y(3) , G(3))$, with the same uniform bound for the zipping lengths extended now to $\Psi (G(2)) = \Phi (G(2)) \subset \Psi (G(3))$, {\it a.s.o.}, with $n \to \infty$, leading to a gigantic $Y(\infty)$. But of course, one needs to show that this kind of thing converges.  It turns out that the techniques from \cite{Po-Ta2} are also beautifully adaptable for this task of rendering convergent the process leading to uniformly bounded zipping length.

\medskip

So far, our representations have been 3-dimensional, but the next step is to go to a very dense 2-dimensional spine of $Y(\infty)$ and then restrict $G(\infty)$ to it. This will lead to Lemma~5 below, which in turn will open the road for the definition of the $S_u (\tilde M^3 (\Gamma))$, already mentioned. The reason for disguising $\Gamma$ as a singular 3-manifold $M^3 (\Gamma)$, was exactly to get all these things straight. But some more terminology is needed before we can go on. For any map $A \overset{h}{\longrightarrow} B$, we will denote by $M_2 (h) \subset A$ the double points set, {\it i.e.} the subset of those $x \in A$ for which ${\rm card} \, (h^{-1} h (x)) > 1$. But then we will also introduce $M^2 (h) \subset A \times A$, which is the set of pairs $(x,y)$, $x \ne y$, {\it s.t.} $h(x) = h(y)$. For instance, $\Phi (f) = M^2 (f) \, \cup \, \mbox{Diag} (X)$. The technology for the next lemma, comes again from the paper \cite{Po-Ta2}.

\bigskip

\noindent {\bf Lemma 5.} {\it There is a $2$-dimensional representation, possessing all the desirable features} 1), 2), 3) {\it from Theorem~$4$, call it now
$$
\xymatrix{
X^2 \ar[rr]^{f}   &&\tilde M^3 (\Gamma)   \\ &
} \leqno (12)
$$
\vglue -8mm
\noindent and which is such that the following things happen too.}

\medskip

4) {\it (The second finiteness condition.) For any tight compact transversal $\Lambda$ to $M_2 (f) \subset X^2$ we have}
$$
\mbox{card} \, (\lim \, (\Lambda \cap M_2 (f))) < \infty \, . \leqno (13)
$$

\medskip

5) {\it The closed subset}
$$
\mbox{LIM} \, M_2 (f) \underset{\rm def}{=} \ \bigcup_{\Lambda} \ \lim \, (\Lambda \cap M_2 (f)) \subset X^2 \leqno (14)
$$
{\it is a locally finite graph and $f \, {\rm LIM} \, M_2 (f) \subset f X^2$ is also a {\ibf closed} subset.}

\medskip

6) {\it Let $\Lambda^*$ run over all tight transversals to ${\rm LIM} \, M_2 (f)$. Then we have}
$$
\bigcup_{\Lambda^*} \ (\Lambda^* \cap M_2 (f)) = M_2 (f) \, .
$$

\bigskip

One should notice that ${\rm LIM} \, M_2 (f) = \emptyset$ is equivalent to $M_2 (f) \subset X^2$ being a closed subset. It so happens that, if this is the case, then it is relatively easy to prove Theorem~1 for the corresponding group $\Gamma$. So, we may as well consider the situation ${\rm LIM} \, M_2 (f) \ne \emptyset$, the only one we will deal with in this paper. But then, see the next section too.

\medskip

Now, once $M_2 (f)$ is NOT a closed subset, the (13) is clearly the next best option. But then, in \cite{Po-Ta3} it is shown that if, forgetting about groups and group actions, we play this same game, in the most economical manner for the Whitehead manifold ${\rm Wh}^3$, instead of $\tilde M^3 (\Gamma)$, then generically the set $\lim (\Lambda \cap M_2 (f))$ becomes a Cantor set, which comes naturally equipped with a feedback loop, turning out to be directly related to the one which generates the Julia sets in the complex dynamics of quadratic maps.

\medskip

Continuing with our list of comments, in the context of 6) in our lemma above, there is a transversal holonomy for ${\rm LIM} \, M_2 (f)$, which is quite non trivial. Life would be easier without that.

\medskip

The $X^2$ has only mortal singularities, while the $fX^2$ only has immortal ones. [In terms of \cite{Po3}, \cite{Ga}, the singularities of $X^2$ (actually of $f$) are all ``undrawable singularities''. At the source, the same is true for the immortal ones too.]

\medskip

The representation space $X^2$ in Lemma~5 is locally finite but, as soon as ${\rm LIM} \, M_2 (f) \ne \emptyset$ (which, as the reader will retain is the main source of headaches in the present paper), the $fX^2 \subset \tilde M^3 (\Gamma)$ is {\ibf not}.

\medskip

At this point, we would like to introduce canonical smooth high-dimensional thickenings for this $fX^2$. Here, even if we temporarily forget about $f \, {\rm LIM} \, M_2 (f)$, there are still the immortal singularities (the only kind which $f X^2$ possesses), to be dealt with; this requires a relatively indirect procedure. One starts with a 4-dimensional smooth thickening (remember that, provisionally, we make as if $f \, {\rm LIM} \, M_2 (f) = \emptyset$). This 4$^{\rm d}$ thickening is not uniquely defined, it depends on a desingularization (see \cite{Ga}). Next one takes the product with $B^m$, $m$ large, and this washes out the desingularization-dependence.

\medskip

In order to deal with ${\rm LIM}Ê\, M_2 (f) \ne \emptyset$, the procedure sketched above has to be supplemented with appropriate {\ibf punctures}, by which we mean here pieces of ``boundary'' of the prospective thickened object which are removed, or sent to infinity. This way we get a more or less canonical {\ibf smooth} high-dimensional thickening for $fX^2$, which we will call $S_u (\tilde M^3 (\Gamma))$. Retain that the definition of $S_u (\tilde M^3 (\Gamma))$ has to include at least a first batch of punctures, just so as to get a smooth object. By now we can also state the

\bigskip

\noindent {\bf Main Lemma 6.} {\it $S_u (\tilde M^3 (\Gamma))$ is} GSC.

\bigskip

For this lemma to be true a second batch of punctures will be necessary. But then, we will also want our $S_u (\tilde M^3 (\Gamma))$ to be such that we should have the implication (2). This will turn out to put a very strict {\ibf ``Stallings barrier''} on how much punctures the definition of $S_u (\tilde M^3 (\Gamma))$ can include, at all. We will have more to say concerning this in Section~4, where some hints concerning the proof of the Main Lemma will also be offered.

\section{A very rough classification for the set of all groups}

The present section, except for some historical hints, corresponds to matters which are even much more ``work in progress'' than Theorem~1, and which will certainly need to be much more worked out, in their details later on. So, with all these cautionary remarks and, at least as a first approximation, let us say that a group $\Gamma$ is {\ibf easy} if it is possible to find for it some 2-dimensional representation with {\ibf closed} $M_2 (f)$. We certainly do not mean here some equivariant representation {\it \`a la} Lemma~5 which, most likely, will have ${\rm LIM} \, M_2 (f) \ne \emptyset$. We do not ask for anything beyond GSC and $\Psi = \Phi$. On the other hand, when there is {\ibf no} representation for $\Gamma$ with a closed $M_2 (f)$, {\it i.e.} if for any 2$^{\rm d}$ representation we have ${\rm LIM} \, M_2 (f) \ne \emptyset$, then we will say that the group $\Gamma$ is {\ibf difficult}.

\medskip

Do not give any connotations to these notions of easy group {\it versus} difficult group, beyond the tentative technical definitions given here.

\medskip

Anyway, for various reasons, we prefer a 3-dimensional alternative definition. All the 3-dimensional representations $X \overset{f}{\longrightarrow} \tilde M^3 (\Gamma)$ are such that $X$ is a union of {\ibf ``fundamental domains''}, pieces on which $f$ injects and which have sizes uniformly bound, both from above and from below. With this, $\Gamma$ will be said now to be easy if for any compact $K \subset \tilde M^3 (\Gamma)$ there is a representation $X \overset{f}{\longrightarrow} \tilde M^3 (\Gamma)$, (possibly depending on $K$), such that only finitely many fundamental domains $\Delta \subset X$ are such that $K \cap f\Delta \ne \emptyset$. There are clearly two distinct notions here, the one just stated and then also the stronger one where a same $(X^2 , f)$ is good for all $K$'s. This last one should certainly be equivalent to the 2-dimensional definition which was given first. Anyway, the general idea here is that the easy groups are those which manage to avoid the {\ibf Whitehead nightmare} which is explained below, and for which we also refer to \cite{Po4}. It was said earlier that we have a not too difficult implication $\{ \Gamma$ is easy$\} \Longrightarrow \{ \Gamma$ is QSF$\}$. I believe this holds even with the $K$-dependent version of ``easy'', but I have not checked this fact.

\medskip

So, the difficult groups are defined now to be those for which {\ibf any} representation $(X,f)$ exhibits the following Whitehead nightmare
$$
\mbox{For any compact $K \subset \tilde M^3 (\Gamma)$ there are INFINITELY many} \leqno (15)
$$
$$
\mbox{fundamental domains $\Delta \subset X$ {\it s.t.} $K \cap f\Delta \ne \emptyset$.}
$$
The Whitehead nightmare above is closely related to the kind of processes {\it via} which the Whitehead manifold ${\rm Wh}^3$ itself or, even more seriously, the Casson Handles are constructed; this is where the name comes from, to begin with.

\medskip

Here is the story behind these notions. Years ago, various people like Andrew Casson, myself, and others, have written papers (of which \cite{Po5}, \cite{Po6}, \cite{G-S},$\ldots$ are only a sample), with the following general gist. It was shown that if for a closed 3-manifold $M^3$, the $\pi_1 \, M^3$ has some kind of nice geometrical features, then $\pi_1^{\infty} \tilde M^3 = 0$; the list of nice geometrical features in question, includes Gromov hyperbolic (or more generally almost convex), automatic (or more generally combable), {\it a.s.o.} The papers just mentioned have certainly been superseded by Perelman's proof of the geometrization conjecture, but it is still instructive to take a look at them from the present vantage point: with hindsight, what they actually did, was to show that, under their respective geometric assumptions, the $\pi_1 \, M^3$ way easy, hence QSF, hence $\pi_1^{\infty} = 0$.

\smallskip

As a small aside, in the context of those old papers mentioned above, both Casson and myself we have developed various group theoretical concepts, out of which Brick and Mihalik eventually abstracted the notion of QSF. For instance, I considered ``Dehn exhaustibility'' which comes with something looking, superficially, like the (QSF) from the statement of Theorem~1 above, except that now both $K$ and $\tilde P$ are smooth and, more seriously so, $F$ is an {\ibf immersion}. Contrary to the QSF itself, the Dehn exhaustibility fails to be presentation independent, but I have still found it very useful now, as an ingredient for the proof of the implication (2). Incidentally also, when groups are concerned, Daniele Otera~\cite{O} has proved that QSF and Dehn exhaustibility are equivalent, in the same kind of weak sense in which he and Funar have proved that QSF and GSC are equivalent.

\medskip

Concerning these same old papers as above, if one forgets about three dimensions and about $\pi_1^{\infty}$ (which I believe to be essentially a red herring in these matters), what they actually prove too, between the lines, is that any $\Gamma$ which satisfies just the geometrical conditions which those papers impose, is in fact easy.

\medskip

So, what next? For a long time I have tried unsuccessfully, to prove that any $\pi_1 \, M^3$ is easy. Today I believe that this is doable, provided one makes use of the geometrization of 3-manifolds in its full glory, {\it i.e.} if one makes full use of Perelman's work. But then, rather recently, I have started looking at these things from a different angle. What the argument for Theorem~1 does, essentially, is to show that even if $\Gamma$ is difficult, it still is QSF. Then, I convinced myself that the argument in question can be twisted around so as to prove a much stronger result, which I only state here as a conjecture, since a lot of details are still to be fully worked out. Here it is:

\bigskip

\noindent {\bf Conjecture 7.} {\it For any $\Gamma$, there is a representation $X \overset{f}{\longrightarrow} \tilde M^3 (\Gamma)$ such that, for any fundamental domain $\delta \subset \tilde M^3 (\Gamma)$, there is a UNIFORMLY BOUNDED number of fundamental domains $\Delta \subset X$ such that $f \Delta \cap \delta \ne \emptyset$.

\medskip

In other words, {\ibf all groups are easy} (in quite a strong sense).}

\bigskip

This certainly implies Theorem~1 and also, according to what I have already said above, it should not really be a new thing for $\Gamma = \pi_1 \, M^3$. But in order to get to that, one has to dig deeper into the details of the Thurston Geometrization, than one needs to do for $\pi_1 \, M^3 \in$ QSF, or at least so I think.

\medskip

Be these things as they may, we leave it now at that and move to an even much more speculative issue.

\medskip

Aperiodic tilings (Penrose tilings, quasi-crystals,$\ldots$) have already received a lot of attention. They are normally considered in Euclidean spaces, {\it i.e.} in the universal covering spaces of tori. But at least in some appropriate conditions it should not be impossible to make sense of them on $\tilde M^n$'s, for instance by projecting down periodic tilings of $\tilde M^n \times R^k$ in an appropriately irrational manner.

\bigskip

\noindent {\bf Question.} Is there some aperiodic tiling on some $\tilde M^n$ such that the $(n-1)$-skeleton is {\ibf not} QSF?

\bigskip

We end the present, more speculative sections here.

\section{Some hints about the proof of the Main Lemma}

We will go back now to those ``punctures'' which have been already considered at the very end of Section~2. So, let $Y$ be some low-dimensional object, like $fX^2$ or some 3-dimensional thickening of it. When it is a matter of punctures, these will be put into effect directly at the level of $Y$ (most usually by an appropriate infinite sequence of Whitehead dilatations), before going to higher dimensions, so that the high dimensional thickening should be {\ibf transversally compact} with respect to the $Y$ to be thickened. The reason for this requirement is that when we will want to prove (2), then arguments like in \cite{Po2} will be used (among others), and these ask for transversal compactness.  For instance, if $V$ is a low-dimensional non-compact manifold, then
$$
\mbox{$V \times B^m$ is transversally compact, while $V \times \{ B^m$ with some boundary}
\leqno (16)
$$
$$
\mbox{punctures$\}$, or even worst $V \times R^m$, is not.}
$$

\medskip

This is our so-called ``Stallings barrier'', putting a limit to how much punctures we are allowed to use. The name is referring here to a corollary of the celebrated Engulfing Theorem of John Stallings, saying that under appropriate dimensions, if $V$ is an {\ibf open} contractible manifold, then $V \times R^p$ is a Euclidean space. There are, of course, also infinitely more simple-minded facts which give GSC when multiplying with $R^p$ or $B^p$.

\medskip

We take now a closer look at the $n$-dimensional smooth manifold $S_u (\tilde M^3 (\Gamma))$ occuring in the Main Lemma. Here $n=m+4$ with high $m$. Essentially, we get our $S_u (\tilde M^3 (\Gamma))$ starting from an initial GSC low-dimensional object, like $X^2$, by performing first a gigantic quotient-space operation, namely our zipping, and finally thickening things into something of dimension $n$. But the idea which we will explore now, is to construct another smooth $n$-dimensional manifold, essentially starting from an already $n$-dimensional GSC smooth manifold, and then use this time a gigantic collection of additions and inclusion maps, which should somehow mimic the zipping process. We will refer to this kind of thing as the {\ibf geometric realization of the zipping}. The additions which are allowed in this kind of realization are Whitehead dilatations, or adding handles of index $\lambda \geq 2$. The final product of the geometric realization, will be another manifold of the same dimension $n$ as $S_u (\tilde M^3 (\Gamma))$, which we will call $S_b (\tilde M^3 (\Gamma))$ and which, {\it a priori} could be quite different from the $S_u (\tilde M^3 (\Gamma))$. To be more precise about this, in a world with ${\rm LIM} \, M_2 (f) = \emptyset$ we would quite trivially have $S_u = S_b$ but, in our real world with ${\rm LIM} \, M_2 (f) \ne \emptyset$, there is absolutely no {\it a priori} reason why this should be so. Incidentally, the subscripts ``$u$'' and ``$b$'', refer respectively to ``usual'' and ``bizarre''.

\medskip

Of course, we will want to compare $S_b (\tilde M^3 (\Gamma))$ and $S_u (\tilde M^3 (\Gamma))$. In order to give an idea of what is at stake here, we will look into the simplest possible local situation with ${\rm LIM} \, M_2 (f)$ present. Ignoring now the immortal singularities for the sake of the exposition, we consider a small smooth chart $U = R^3 = (x,y,z)$ of $\tilde M^3 (\Gamma)$, inside which live $\infty + 1$ planes, namely $W = (z=0)$ and the $V_n = (x=x_n)$, where $x_1 < x_2 < x_3 < \ldots$ with $\lim x_n = x_{\infty}$. Our local model for $X^2 \overset{f}{\longrightarrow} \tilde M^3 (\Gamma)$ is here $f^{-1} \, U = W + \underset{1}{\overset{\infty}{\sum}} \ V_n \subset X^2$ with $f \mid \Bigl( W + \underset{1}{\overset{\infty}{\sum}} \ V_n \Bigl)$ being the obvious map. We find here that the line $(x = x_{\infty}, z=0) \subset W$ is in ${\rm LIM} \, M_2 (f)$ and the situation is sufficiently simple so that we do not need to distinguish here between ${\rm LIM} \, M_2 (f)$ and $f \, {\rm LIM} \, M_2 (f)$.

\medskip

Next, we pick up a sequence of positive numbers converging very fast to zero $\varepsilon > \varepsilon_1 > \varepsilon_2 > \ldots$ and, with this, on the road to the $S_u (\tilde M^3 (\Gamma))$ from the {\ibf Main Lemma}, we will start by replacing the $f W \cup \underset{1}{\overset{\infty}{\sum}} \ V_n \subset f X^2$, with the following 3-dimensional non-compact 3-manifold with boundary
$$
M \underset{\rm def}{=} [ W \times (-\varepsilon \leq z \leq \varepsilon) - {\rm LIM} \, M_2 (f) \times \{ z = \pm \, \varepsilon \}] \ \cup  \leqno (17)
$$
$$
\cup \ \sum_{1}^{\infty} V_n \times (x_n - \varepsilon_n \leq x \leq x_n + \varepsilon_n ) \, .
$$
In such a formula, notations like ``$W \times (-\varepsilon \leq z \leq \varepsilon)$'' should be read ``$W$ thickened into $-\varepsilon \leq z \leq \varepsilon$''. Here ${\rm LIM} \, M_2 (f) \times \{ \pm \, \varepsilon \}$ is a typical puncture, necessary to make our $M$ be a smooth manifold. For expository purposes, we will pretend now that $n=4$ and then $M \times [0 \leq t \leq 1]$ is a local piece of $S_u (\tilde M^3 (\Gamma))$. Now, in an ideal world (but not in ours!), the geometrical realization of the zipping process {\it via} the inclusion maps (some of which will correspond to the Whitehead dilatations which are necessary for the punctures), which are demanded by $S_b (\tilde M^3 (\Gamma))$, should be something like this. We start with the obviously GSC $n$-dimensional thickening of $X^2$, call it $\Theta^n (X^2)$; but remember that for us, here, $n=4$. Our local model should live now inside $R^4 = (x,y,z,t)$, and we will try to locate it there conveniently for the geometric realization of the zipping. We will show how we would like to achieve this for a generic section $y = $ constant.

\medskip

For reasons to become soon clear, we will replace the normal section $y = $ const corresponding to $W$ and which should be
$$
N_y = [ - \infty < x < \infty \, , \ y = {\rm const} \, , \ - \varepsilon \leq z \leq \varepsilon \, , \ 0 \leq t \leq 1] \leqno (18)
$$
$$
- \, (x = x_{\infty} \, , \ y = {\rm const} \, , \ z = \pm \, \varepsilon \, , \ 0 \leq t \leq 1) \, ,
$$
by the smaller $N_y \, - \, \overset{\infty}{\underset{1}{\sum}} \ {\rm DITCH} \, (n)_y$, which is defined as follows. The ${\rm DITCH} \, (n)_y$ is a thin column of height $-\varepsilon \leq z \leq \varepsilon$ and of $(x,t)$-width $4 \, \varepsilon_n$, which is concentrated around the arc
$$
(x = x_n \, , \ y = {\rm const} \, , \ - \varepsilon \leq z \leq \varepsilon \, , \ t=1) \, .
$$
This thin indentation inside $N_y$ is such that, with our fixed $y = {\rm const}$ being understood here, we should have
$$
\lim_{n = \infty} \, {\rm DITCH} \, (n)_y = (x = x_{\infty} \, , \ -\varepsilon \leq z \leq \varepsilon \, , \ t=1) \, . \leqno (19)
$$
Notice that, in the RHS of (19) it is exactly the $z = \pm \, \varepsilon$ which corresponds to punctures.

\medskip

Continuing to work here with a fixed, generic $y$, out of the normal $y$-slice corresponding to $V_n$, namely $(x_n - \varepsilon_n \leq x \leq x_n + \varepsilon_n \, , \ - \infty < z < \infty \, , \ 0 \leq t \leq 1)$, we will keep only a much thinner, isotopically equivalent version, namely the following
$$
(x_n - \varepsilon_n \leq x \leq x_n + \varepsilon_n \, , \ - \infty < z < \infty \, , \ 1 - \varepsilon_n \leq t \leq 1) \, . \leqno (20)
$$
This (20) has the virtue that it can fit now inside the corresponding ${\rm DITCH} \, (n)_y$, without touching at all the $N_y - \{{\rm DITCHES}\}$.

\medskip

What has been carefully described here, when all $y$'s are being taken into account, is a very precise way of separating the $\infty + 1$ branches of (the thickened) (17), at the level of $R^4$, taking full advantage of the additional dimensions ({\it i.e.} the factor $[0 \leq t \leq 1]$ in our specific case). With some work, this kind of thing can be done consistently for the whole global $fX^2$. The net result is an isotopically equivalent new model for $\Theta^n (X^2)$, which invites us to try the following naive approach for the geometric realization of the zipping. Imitating the successive folding maps of the actual zipping, fill up all the empty room left inside the ditches, by using only Whitehead dilatations and additions of handles of index $\lambda > 1$, until one has reconstructed completely the $S_u (\tilde M^3 (\Gamma))$. Formally there is no obstruction here and then also what at a single $y = $ const may look like a handle of index one, becomes ``only index $\geq 2$'', once the full global zipping is taken into account. But there {\ibf is} actually a big problem with this naive approach, {\it via} which one can certainly reconstruct $S_u (\tilde M^3 (\Gamma))$ as a set, but with the {\ibf wrong topology}, as it turns out. I will give an exact idea now of how far one can actually go, proceeding in this naive way. In \cite{Po-Ta4} we have tried to play the naive game to its bitter end, and the next Proposition~8, given here for purely pedagogical reasons, is the kind of thing one gets (and certainly nothing better than it).

\bigskip

\noindent {\bf Proposition 8.} {\it Let $V^3$ be any open simply-connected $3$-manifold, and let also $m \in Z_+$ be high enough. There is then an infinite collection of smooth $(m+3)$-dimensional manifolds, all of them non-compact, with very large boundary, connected by a sequence of smooth embeddings
$$
X_1 \subset X_2 \subset X_3 \subset \ldots \leqno (21)
$$
such that}

\medskip

1) {\it $X_1$ is {\rm GSC} and each of the inclusions in {\rm (21)} is either an elementary Whitehead dilatation or the addition of a handle of index $\lambda > 1$.}

\medskip

2) {\it When one considers the union of the objects in {\rm (21)}, endowed with the {\ibf weak topology}, and there is no other, reasonable one which is usable here, call this new space $\varinjlim \, X_i$, then there is a continuous bijection}
$$
\xymatrix{
\varinjlim \, X_i  \ar[rr]^{\psi}   &&V^3 \times B^m \, .   \\ &
} \leqno (22)
$$
\vglue -5mm
The reader is reminded that in the weak topology, a set $F \subset \varinjlim \, X_i$ is closed iff all the $F \cap X_i$ are closed. Also, the inverse of $\psi$ is not continuous here; would it be, this would certainly contradict \cite{Po2}, since $V^3$ may well be ${\rm Wh}^3$. This, {\it via} Brouwer, also means that $\varinjlim \, X_i$ cannot be a manifold (which would automatically be then GSC). Also, exactly for the same reasons why $\varinjlim \, X_i$ is not a manifold, it is not a metrizable space either. So, here we meet a new barrier, which I will call the {\ibf non metrizability barrier} and, when we will really realize geometrically the zipping, we better stay on the good side of it.

\medskip

One of the many problems in this paper is that the Stallings barrier and the non metrizability barrier somehow play against each other, and it is instructive to see this in a very simple instance. At the root of the non metrizability are, as it turns out, things like (19); {\it a priori} this might, conceivably, be taken care of by letting all of $(x = x_{\infty} \, , \ - \varepsilon \leq z \leq \varepsilon \, , \ t=1)$ be punctures, not just the $z = \pm \, \varepsilon$ part. But then we would also be on the wrong side of the Stallings barrier. This kind of conflict is quite typical.

\medskip

So far, we have presented the disease and the rest of the section gives, essentially, the cure. In a nutshell here is what we will do. We start by drilling a lot of {\ibf Holes}, consistently, both at the (thickened) levels of $X^2$ and $fX^2$. Working now only with objects with Holes, we will be able to fill in the empty space left inside the ${\rm Ditch} \, (n)$ {\ibf only} for $1-\varepsilon_n \leq z \leq 1$ where, remember $\underset{n = \infty}{\rm lim} \, \varepsilon_n = 0$. This replaces the trouble-making (19) by the following
$$
\lim_{n = \infty} \{\mbox{{\it truncated}} \ {\rm Ditch} \, (n)_y \} = (x = x_{\infty} \, , \ z = \varepsilon \, , \ t = 1) \, , \leqno (23)
$$
which is now on the good side, both of the Stallings barrier and of the non metrizability barrier. But, after this {\ibf partial ditch-filling process}, we have to go back to the Holes and put back the corresponding deleted material. This far from trivial issue will be discussed later.

\medskip

But before really going on, I will make a parenthetical comment which some readers may find useful. There are many differences between the present work (of which this paper is only a summary, the long complete version exists too, in hand-written form, waiting to be typed), and my ill-fated, by now dead $\pi_1^{\infty} \tilde M^3 = 0$ attempt [Pr\'epublication Orsay 2000-20 and 2001-57]. Of course, a number of ideas from there found their way here too. In the dead papers I was also trying to mimick the zipping by inclusions; but there, this was done by a system of ``gutters'' added in $2^{\rm d}$ or $3^{\rm d}$, before any thickening into high dimensions. Those gutters came with fatal flaws, which opened a whole Pandora's box of troubles. These turned out to be totally unmanageable, short of some input of new ideas. And, by the time a first whiff of such ideas started popping up, Perelman's announcement was out too. So, I dropped the whole thing, turning to more urgent tasks.

\medskip

It is only a number of years later that this present work grew out of the shambles. By contrast with the low-dimensional gutters from the old discarded paper, the present ditches take full advantages of the {\ibf additional} dimensions; I use here the word ``additional'', by opposition to the mere high dimensions. The spectre of non metrizability which came with the ditches, asked then for Holes the compensating curves of which are dragged all over the place by the inverse of the zipping flow, far from their normal location, {\it a.s.o.} End of practice.

\medskip

With the Holes in the picture, the $\Theta^n (X^2)$, $S_u (\tilde M^3 (\Gamma))$, will be replaced by the highly non-simply-connected smooth $n$-manifolds with non-empty boundary $\Theta^n (X^2) - H$ and $S_u (\tilde M^3 (\Gamma)) - H$. Here the ``$-H$'' stands for ``with the Holes having been drilled'' or, in a more precise language, as it will turn out, with the 2-handles which correspond to the Holes in question, deleted.

\medskip

The manifold $S_u (\tilde M^3 (\Gamma)) - H$ comes naturally equipped with a PROPER framed link
$$
\xymatrix{
\overset{\infty}{\underset{1}{\sum}} \, C_n \ar[rr]^{\!\!\!\!\!\!\!\!\!\!\!\!\!\!\!\alpha}   &&\partial \, (S_u (\tilde M^3 (\Gamma)) - H) \, ,    \\ &
} \leqno (24)
$$
\vglue -8mm
\noindent where ``$C$'' stands for ``curve''. This framed link is such that, when one adds the corresponding 2-handles to $S_u (\tilde M^3 (\Gamma)) - H$, then one gets back the $S_u (\tilde M^3 (\Gamma))$.

\medskip

This was the $S_u$ story with holes, which is quite simple-minded, and we move now on to $S_b$. One uses now the $\Theta^n (X^2) - H$, {\it i.e.} the thickened $X^2$, with Holes, as a starting point for the geometric realization of the zipping process, rather than starting from the $\Theta^n (X^2)$ itself. The Holes allow us to make use of a partial ditch filling, {\it i.e.} to use now the truncation $1-\varepsilon_n \leq z \leq 1$, which was already mentioned before. This has the virtue of putting us on the good side of all the various barriers which we have to respect. The end-product of this process is another smooth $n$-dimensional manifold, which we call $S_b (\tilde M^3 (\Gamma)) - H$. This comes with a relatively easy diffeomorphism
$$
\xymatrix{
S_b (\tilde M^3 (\Gamma)) - H \ar[rr]_{\approx}^{\eta}   &&S_u (\tilde M^3 (\Gamma)) - H \, .   \\ &
} \leqno (25)
$$
\vglue -8mm
\noindent Notice that only ``$S_b (\tilde M^3 (\Gamma)) - H$'' is defined, so far, and not yet the full $S_b (\tilde M^3 (\Gamma))$ itself. Anyway here comes the following fact, which is far from being trivial.

\bigskip

\noindent {\bf Lemma 9.} {\it There is a second PROPER framed link
$$
\xymatrix{
\underset{1}{\overset{\infty}{\sum}} \, C_n \ar[rr]^{\!\!\!\!\!\!\!\!\!\!\!\!\!\!\!\beta}   &&\partial \, (S_b (\tilde M^3 (\Gamma)) - H)   \\ &
} \leqno (26)
$$
\vglue -8mm
\noindent which has the following two features.}

\medskip

1) {\it The following diagram is commutative, {\ibf up to a homotopy}}
$$
\xymatrix{
S_b (\tilde M^3 (\Gamma)) -H \ar[rr]_{\eta}   &&S_u (\tilde M^3 (\Gamma)) - H   \\ 
&\overset{\infty}{\underset{1}{\sum}} \, C_n \ar[ul]^{\beta} \ar[ur]_{\alpha}
} \leqno (27)
$$
{\it The homotopy above, which is {\ibf not} claimed to be PROPER, is compatible with the framings of $\alpha$ and $\beta$.}

\medskip

2) {\it We {\rm define} now the smooth $n$-dimensional manifold
$$
S_b (\tilde M^3 (\Gamma)) = (S_b (\tilde M^3 (\Gamma)) - H) + \{\mbox{\rm the $2$-handles which are defined by} \leqno (28)
$$
$$
\mbox{\rm the framed link $\beta$ {\rm (26)}}\} \, .
$$
This manifold $S_b (\tilde M^3 (\Gamma))$ is} GSC.

\bigskip

Without even trying to prove anything, let us just discuss here some of the issues which are involved in this last lemma.

\medskip

In order to be able to discuss the (26), let us look at a second toy-model, the next to appear, in increasing order of difficulty, after the one already discussed, when formulae (17) to (20) have been considered. We keep now the same $\underset{1}{\overset{\infty}{\sum}} \, V_n$, and just replace the former $W$ by $W_1 = (y=0) \cup (z=0)$. The $M$ from (17) becomes now the following non-compact 3-manifold with boundary
$$
M_1 = [(-\varepsilon \leq y \leq \varepsilon) \cup (-\varepsilon \leq z \leq \varepsilon) - \{\mbox{the present contribution of} \leqno (29)
$$
$$
{\rm LIM} \, M_2 (f) \}] \cup \sum_{1}^{\infty} \, V_n \times (x_n - \varepsilon_n \leq x \leq x + \varepsilon_n ) \, ;
$$
the reader should not find it hard to make explicit the contribution of ${\rm LIM} \, M_2 (f)$ here. Also, because we have considered only $(y=0) \, \cup \, (z=0)$ and not the slightly more complicated disjoint union $(y=0) + (z=0)$, which comes with triple points, there is still no difference, so far, between ${\rm LIM} \, M_2 (f)$ and $f \, {\rm LIM} \, M_2 (f)$.

\medskip

When we have discussed the previous local model, then the
$$
{\rm DITCH} \, (n) \underset{\rm def}{=} \ \bigcup_y \ {\rm DITCH} \, (n)_y
$$
was concentrated in the neighbourhood of the rectangle
$$
(x = x_n \, , \ -\infty < y < \infty \, , \ - \varepsilon \leq z \leq \varepsilon \, , \ t=1 ) \, .
$$
Similarly, the present ${\rm DITCH} \, (n)$ will be concentrated in a neighbourhood of the $2$-dimensional infinite cross
$$
(x = x_n \, , \ (-\varepsilon \leq y \leq \varepsilon) \cup (-\varepsilon \leq z \leq \varepsilon) \, , \ t=1) \, .
$$
It is only the $V_n$'s which see Holes. Specifically, $V_n - H$ is a very thin neighbourhood of the $1^{\rm d}$ cross $(y = + \varepsilon) \cup (z = +\varepsilon)$, living at $x=x_n$, for some fixed $t$, and fitting inside ${\rm DITCH} \, (n)$, without touching anything at the level of the $\{ W_1$ thickened in the high dimension, and with the DITCH deleted$\}$. But it does touch to four Holes, corresponding to the four corners. The action takes place in the neighbourhood of $t=1$, making again full use of the additional dimensions, supplementary to those of $M_1$ (29).

\medskip

With this set-up, when we try to give a ``normal'' definition for the link $\beta$ in (26), then we encounter the following kind of difficulty, and things only become worse when triple points of $f$ are present too. Our normal definition of $\beta \, C_n$ (where $C_n$ is here the generic boundary of one of our four Holes), is bound to make use of arcs like $I_n = (x=x_n \, , \ y = +\varepsilon \, , \ -\varepsilon \leq z \leq \varepsilon \, , \ t = {\rm const})$, or $I_n = (x = x_n \, , \ - \varepsilon \leq y \leq \varepsilon \, , \ z = + \varepsilon \, , \ t = {\rm const})$ which, in whatever $S_b (\tilde M^3 (\Gamma)) - H$ may turn out to be, accumulate at finite distance. So, the ``normal'' definition of $\beta$ fails to be PROPER, and here is what we will do about this. Our arcs $I_n$ come naturally with (not completely uniquely defined) double points living in $\Psi (f) = \Phi (f)$, attached to them and these have zipping paths, like in Theorem~4. The idea is then to push the arcs $I_n$ back, along (the inverses of) the zipping paths, all the way to the singularities of $f$, and these do not accumulate at finite distance. In more precise terms, the arcs via which we replace the $I_n$'s inside $\beta \, C_n$ come closer and closer to $f \, {\rm LIM} \, M_2 (f)$ as $n \to \infty$, and this last object lives at infinity. This is the {\it correct} definition of $\beta$ in (26). The point is that, now $\beta$ is PROPER, as claimed in our lemma. But then, there is also a relatively high price to pay for this. With the roundabout way to define $\beta$ which we have outlined above, the point 2) in our lemma, {\it i.e.} the basic property $S_b \in$ GSC, is no longer the easy obvious fact which it normally should be; it has actually become non-trivial. Here is where the difficulty sits. To begin with, our $X^2$ which certainly is GSC, by construction, houses two completely distinct, not everywhere well-defined flow lines, namely the zipping flow lines and the collapsing flow lines stemming from the easy id $+$ nilpotent geometric intersection matrix of $X^2$. By itself, each of these two systems of flow lines is quite simple-minded and controlled, but not so the combined system, which can exhibit, generally speaking, closed oriented loops.

\medskip

These {\ibf bad cycles} are in the way for the GSC property of $S_b$; they introduce unwanted terms in the corresponding geometric intersection matrix. An infinite machinery is required for handling this new problem: the bad cycles have to be, carefully, pushed all the way to infinity, out of the way. Let me finish with Lemma~9 by adding now the following item. Punctures are normally achieved by infinite sequences of dilatations, but if we locate the $\beta \, \underset{1}{\overset{\infty}{\sum}} \, C_n$ over the regions created by them, this again may introduce unwanted terms inside the geometric intersection matrix of $S_b$, making havoc of the GSC feature. In other words, we have now additional restrictions concerning the punctures, going beyond what the Stallings barrier would normally tolerate. In more practical terms, what this means that the use of punctures is quite drastically limited, by the operations which make our $\beta$ (26) PROPER. This is as much as I will say about the Lemma~9, as such.

\medskip

But now, imagine for a minute that, in its context, we would also know that (27) commutes up to PROPER homotopy. In that hypothetical case, in view of the high dimensions which are involved, the (27) would commute up to isotopy too. This would prove then that the manifolds $S_u (\tilde M^3 (\Gamma))$ and $S_b (\tilde M^3 (\Gamma))$ are diffeomorphic. Together with 2) in Lemma 9, we would then have a proof of our desired Main Lemma. But I do not know how to make such a direct approach work, and here comes what I can offer instead.

\medskip

The starting point of the whole $S_u / S_b$ story, was an {\ibf equivariant} representation theorem for $\tilde M^3 (\Gamma)$, namely our Theorem~4 and, from there on, although we have omitted to stress this until now, everything we did was supposed to be equivariant all along: the zipping, the Holes, the (24) $+$ (26), {\it a.s.o.} are all equivariant things. 

\medskip

Also, all this equivariant discussion was happening upstairs, at the level of $\tilde M^3 (\Gamma)$. But, being equivariant it can happily be pushed down to the level of $M^3 (\Gamma) = \tilde M^3 (\Gamma) / \Gamma$. Downstairs too, we have now two, still non-compact manifolds, namely $S_u (M^3 (\Gamma)) \underset{\rm def}{=} S_u (\tilde M^3 (\Gamma)) / \Gamma$ and $S_b (M^3 (\Gamma)) \underset{\rm def}{=} S_b (\tilde M^3 (\Gamma)) / \Gamma$.

\medskip

What these last formulae mean, is also that we have
$$
S_u (M^3 (\Gamma))^{\sim} = S_u (\tilde M^3 (\Gamma)) \quad \mbox{and} \quad S_b (M^3 (\Gamma))^{\sim} = S_b (\tilde M^3 (\Gamma)) \, , \leqno (30)
$$
let us say that $S_u$ and $S_b$ are actually functors of sorts.

\medskip

But before really developping this new line of thought, we will have to go back to the diagram (27) which, remember, commutes up to homotopy. Here, like in the elementary text-books, we would be very happy now to change
$$
\alpha \, C_n \sim \eta \beta \, C_n \quad \mbox{into (something like)} \quad \alpha \, C_n \cdot \eta \beta \, C_n^{-1} \sim 0 \, .
$$
This is less innocent than it may look, since in order to be of any use for us, the infinite system of curves
$$
\Lambda_n \underset{\rm def}{=} \alpha \, C_n \cdot \eta \beta \, C_n^{-1} \subset \partial \, (S_u (\tilde M^3 (\Gamma)) - H) \, , \ n = 1,2, \ldots \leqno (31)
$$
better be PROPER (and, of course, equivariant too). The problem here is the following and, in order not to over complicate our exposition, we look again at our simplest local model. Here, in the most difficult case at least, the curve $\alpha \, C_n$ runs along $(x = x_n \, , \ z = - \varepsilon)$ while $\eta \beta \, C_n$ runs along $(x = x_n \, , \ z = + \varepsilon)$. The most simple-minded procedure for defining (31) would then be to start by joining them along some arc of the form
$$
\lambda_n = (x = x_n \, , \ y = {\rm const} \, , \ -\varepsilon \leq z \leq \varepsilon \, , \ t = {\rm const}) \, . \leqno (32)
$$

But then, for the very same reasons as in our previous discussion of the mutual contradictory effects of the two barriers (Stallings and non-metrizability), this procedure is certainly not PROPER. The cure for this problem is to use, once again, the same trick as for defining a PROPER $\beta$ in (26), namely to push the stupid arc $\lambda_n$ along the (inverse of) the zipping flow, all the way back to the singularities, keeping things all the time close to $f \, {\rm LIM} \, M_2 (f)$, {\it i.e.} close to $(x=x_{\infty} \, , \ z = \pm \varepsilon \, , \ t=1)$, in the beginning at least.

\medskip

The next lemma sums up the net result of all these things.

\bigskip

\noindent {\bf Lemma 10.} {\it The correctly defined system of curves {\rm (31)} has the following features}
\begin{itemize}
\item[0)] {\it It is equivariant (which, by now, does not cost much),}
\item[1)] {\it It is PROPER,}
\item[2)] {\it For each individual $\Lambda_n$, we have a null homotopy
$$
\Lambda_n \sim 0 \quad \mbox{in} \quad S_u (\tilde M^3 (\Gamma)) - H \, .
$$
[Really this is in $\partial (S_u (\tilde M^3 (\Gamma)) - H)$, but we are very cavalier now concerning the distinction between $S_u$ and $\partial S_u$; the thickening dimension is very high, anyway.]}
\item[3)] {\it As a consequence of the bounded zipping length in the Theorem~$4$, our system of curves $\Lambda_n$ has {\ibf uniformly bounded length}.}
\end{itemize}

\bigskip

Lemma~10 has been stated in the context of $\tilde M^3 (\Gamma)$, upstairs. But then we can push it downstairs to $M^3 (\Gamma)$ too, retaining 1), 2), 3) above. So, from now on, we consider the correctly defined system $\Lambda_n$ {\it downstairs}. Notice here that, although $M^3 (\Gamma)$ is of course compact (as a consequence of $\Gamma$ being finitely presented), the $S_u (M^3 (\Gamma)) - H$ and $S_u (M^3 (\Gamma))$ are certainly not. 

\medskip

So, the analogue of 1) from Lemma~10, which reads now
$$
\lim_{n = \infty} \, \Lambda_n = \infty \, , \ \mbox{inside} \ S_u (M^3 (\Gamma))-H \, , \leqno (33)
$$
is quite meaningful. Finally, here is the

\bigskip

\noindent {\bf Lemma 11.} (KEY FACT) {\it The analogue of diagram {\rm (27)} downstairs, at the level of $M^3 (\Gamma)$, commutes now up to PROPER homotopy.}

\bigskip

Before discussing the proof of this key fact, let us notice that it implies that $S_u (M^3 (\Gamma)) \underset{\rm DIFF}{=} S_b (M^3 (\Gamma))$ hence, {\it via} (30), {\it i.e.} by ``functoriality'', we also have
$$
S_u (\tilde M^3 (\Gamma)) \underset{\rm DIFF}{=} S_b (\tilde M^3 (\Gamma)) \, , \leqno (34)
$$
making that $S_u (\tilde M^3 (\Gamma))$ is GSC, as desired; see here 2) in Lemma~9 too.

\medskip

All the rest of the discussion is now downstairs, and we will turn back to Lemma~10. Here, the analogue of 2) is, of course, valid downstairs too, which we express as follows
$$
\mbox{For every $\Lambda_n$ there is a singular disk $D_n^2 \subset S_u (M^3 (\Gamma)) - H$,} \leqno (35)
$$
$$
\mbox{with $\partial D^2 = \Lambda_n$.} 
$$
With this, what we clearly need now for Lemma~11, is something like (35), but with the additional feature that $\underset{n=\infty}{\lim} \, D_n^2 = \infty$, inside $S_u (M^3 (\Gamma)) - H$. In a drastically oversimplified form, here is how we go about this. Assume, by contradiction, that there is a compact set $K \subset \partial (S_u (M^3 (\Gamma)) - H)$ and a subsequence of $\Lambda_1 , \Lambda_2 , \ldots$, which we denote again by exactly the same letters, such that for {\ibf any} corresponding system  of singular disks cobounding it, $D_1^2 , D_2^2 , \ldots$, we should have $K \cap D_n^2 \ne \emptyset$, for all $n$'s.

\medskip

We will  show now that this itself, leads to a contradiction. Because $M^3 (\Gamma)$ is compact ($\Gamma$ being finitely presented), we can {\ibf compactify} $S_u (M^3 (\Gamma)) - H$ by starting with the normal embedding $S_u (M^3 (\Gamma)) - H \subset M^3 (\Gamma) \times B^N$, $N$ large; inside this compact metric space, we take then the closure of $S_u (M^3 (\Gamma)) - H$. For this compactification which we denote by $(S_u (M^3 (\Gamma)) - H)^{\wedge}$ to be nice and useful for us, we have to be quite careful about the exact locations and sizes of the Holes, but the details of this are beyond the present outline. This good compactification is now
$$
(S_u (M^3 (\Gamma)) - H)^{\wedge} = (S_u (M^3 (\Gamma)) - H) \cup E_{\infty} \, , \leqno (36)
$$
where $E_{\infty}$ is the compact space which one has to add at the infinity of $S_u (M^3 (\Gamma))$ $- \, H$, so as to make it compact. It turns out that $E_{\infty}$ is moderately wild, failing to the locally connected, although it has plenty of continuous arcs embedded inside it.

\medskip

Just by metric compactness, we already have
$$
\lim_{n = \infty} {\rm dist} (\Lambda_n , E_{\infty}) = 0 \leqno (37)
$$
and, even better, once we know that the lengths of $\Lambda_n$ are uniformly bounded, there is a subsequence $\Lambda_{j_1} , \Lambda_{j_2} , \Lambda_{j_3} , \ldots$ of $\Lambda_1 , \Lambda_2 , \Lambda_3 , \ldots$ and a continuous curve $\Lambda_{\infty} \subset E_{\infty}$, such that $\Lambda_{j_1} , \Lambda_{j_2} , \ldots$ converges uniformly to $E_{\infty}$. To be pedantically precise about it, we have
$$
{\rm dist} \, (\Lambda_{j_n} , \Lambda_{\infty}) = \varepsilon_n \, , \ \mbox{where} \ \varepsilon_1 > \varepsilon_2 > \ldots > 0 \ \mbox{and} \ \lim_{n = \infty} \varepsilon_n = 0 \, . \leqno (38)
$$
Starting from this data, and injecting also a good amount of precise knowledge concerning the geometry of $S_u (M^3 (\Gamma))$ (knowledge which we actually have to our disposal, in real life), we can construct a region $N = N (\varepsilon_1 , \varepsilon_2 , \ldots) \subset S_u (M^3 (\Gamma)) - H$, which has the following features

\medskip

A) The map $\pi_1 N \longrightarrow \pi_1 (S_u (M^3 (\Gamma)) - H)$ injects;

\medskip

B) The $E_{\infty}$ lives, also, at the infinity of $N$, and there is a retraction
$$
\xymatrix{
N \cup E_{\infty} \ar[rr]^{R}   &&E_{\infty} \, ;   \\ &
}
$$
\vglue -8mm
C) There is an ambient isotopy of $S_u (M^3 (\Gamma)) - H$, which brings all the $\Lambda_{j_1} , \Lambda_{j_2} , \ldots$ inside $N$. After this isotopy, we continue to have (38), or at least something very much like it.

\medskip

The reader may have noticed that, for our $N$ we studiously have avoided the word ``neighbourhood'', using ``region'' instead; we will come back to this.

\medskip

Anyway, it follows from A) above that our $\Lambda_{j_1} , \Lambda_{j_2} , \ldots \subset N$ (see here C)) bound singular disks in $N$. Using B), these disks can be brought very close to $E_{\infty}$, making them disjoined from $K$. In a nutshell, this is the contradiction which proves what we want. We will end up with some comments.

\medskip

To begin with, our $N = N (\varepsilon_1 , \varepsilon_2 , \ldots)$ is by no means a neighbourhood of infinity, it is actually too thin for that, and its complement is certainly not pre-compact.

\medskip

In the same vein, our argument which was very impressionistically sketched above, is certainly not capable of proving things like $\pi_1^{\infty} (S_u (M^3 (\Gamma)) - H) = 0$, which we do {\ibf not} claim, anyway.

\medskip

Finally, it would be very pleasant if we could show that $\Lambda_{\infty}$ bounds a singular disk inside $E_{\infty}$ and deduce then our desired PROPER homotopy for the $\Lambda_{i_n}$'s from this. Unfortunately, $E_{\infty}$ is too wild a set to allow such an argument to work.

\end{document}